\def\IN{\mathaccent "017}
\documentclass[12pt]{amsart}
\usepackage{amsmath,amssymb,amscd}
\newtheorem{theorem}{Theorem}
\newtheorem{proposition}{Proposition}
\newtheorem{lemma}{Lemma}
\newtheorem{corollary}{Corollary}
\begin{document}
\centerline{\Large Lightlike foliations of semi-Riemannian
manifolds\footnote{ \vskip 0.02in MS
classification: 53C12, 53C50.}} \vskip 1.5cm \centerline{\large
Elisabetta Barletta \hspace{0.5cm} Sorin
Dragomir\footnote{Universit\`a della Basilicata, Dipartimento di
Matematica, Contrada Macchia Romana, 85100 Potenza, Italy, e-mail:
{\tt barletta@unibas.it}, {\tt dragomir@unibas.it}}} \vskip 0.5cm
\centerline{\large Krishan L. Duggal\footnote{University of
Windsor, Department of Mathematics and Statistics, Windsor,
Ontario N9B3P4, Canada, http://www.uwindsor.ca/duggal, e-mail:
{\tt yq8@uwindsor.ca}}}

\title{}
\author{}
\begin{abstract} Using screen distributions and lightlike transversal vector bundles
we develop a theory of degenerate foliations of semi-Riemannian
manifolds. We build lightlike foliations of a semi-Riemannian
manifold by suspension of a group homomorphism $\varphi : \pi_1
(B, x_0 ) \to {\rm Isom}(T)$. We compute the basic cohomology
groups of the flow determined by a lightlike Killing vector field
on a complete semi-Riemannian manifold. We prove a lightlike
analog to Rummler's formula and the transversal divergence theorem
of F. Kamber et al., \cite{kn:KTT}.
\end{abstract}
\maketitle

\section{Introduction}
A lightlike foliation $\mathcal F$ of a semi-Riemannian manifold
$M$ is a foliation each of whose leaves is a lightlike submanifold
on $M$, so that the restriction of the ambient metric to the
tangent bundle $T({\mathcal F})$ is degenerate. Therefore ${\rm
Rad} \; T{\mathcal F} := T({\mathcal F}) \cap T({\mathcal F})^\bot
\neq 0$ and one may not develop a satisfactory theory (a geometry
of the second fundamental form of $\mathcal F$ in $M$) by a mere
imitation of the theory of foliations of Riemannian manifolds, cf.
e.g. \cite{kn:Ton}, p. 62-73. Indeed the very basics (existence of
bundle-like metrics for Riemannian foliations, building adapted
connections in the normal bundle, etc.) depend upon the
availability of a natural isomorphism $\sigma : \nu ({\mathcal F})
\approx T({\mathcal F})^\bot$ whose existence follows from the
nondegeneracy of $T({\mathcal F})$. We solve this problem on the
lines of \cite{kn:DuBe} (which deals with the case of a single
lightlike submanifold) by using the technique of screen
distributions and the corresponding lightlike transversal bundles.
Precisely we build a vector bundle $tr(T{\mathcal F}) \to M$
(depending on a choice of complements - the so called {\em screen
distributions} - to ${\rm Rad} \; T{\mathcal F}$ in $T({\mathcal
F})$ and $T({\mathcal F})^\bot$ respectively) playing the role of
$T({\mathcal F})^\bot$ in the theory of foliations with
nondegenerate leaves, i.e. $T(M) = T({\mathcal F}) \oplus
tr(T{\mathcal F})$, whose key property is that its lightlike part
$ltr(T{\mathcal F})$ is {\em not} orthogonal to the radical
distribution. The particular cases we study are those of lightlike
foliations defined by suspension, flows of lightlike Killing
vector fields, and foliations by level sets of lightlike functions
i.e. smooth functions on a semi-Riemannian manifold whose gradient
is null. By a result of Y. Kamishima, \cite{kn:Kam}, a Lorentz
spherical manifold $M$ admits no timelike or lightlike Killing
vector fields (and if $M$ is compact and $3$-dimensional there are
no spacelike Killing vector fields as well). As an application of
our theory, given a complete $3$-dimensional Lorentz manifold we
may weaken the hypothesis in \cite{kn:Kam} by assuming that ${\rm
Isom}(M) = {\rm O}(4,1)$ and that $M$ has the real homology of a
pseudosphere $S^3_1 (r)$ proving however a less precise result:
such $M$ admits no {\em complemented} lightlike Killing vector
fields (cf. Corollary \ref{c:fol1}). More general, given a
lightlike Killing vector field and the corresponding flow
$\mathcal F$ on a complete semi-Riemannian manifold we build a
long exact sequence of cohomology groups
\[ H^k_B ({\mathcal F}) \to H^k (M , {\mathbb R}) \to H^{k-1}_B
({\mathcal F}) \stackrel{\Delta}{\longrightarrow} H^{k+1}_B
({\mathcal F}) \to \cdots
\]
allowing one to compute the basic cohomology of the flow when the
de Rham cohomology of $M$ is known (e.g. when $M \sim S^n_s (r)$
i.e. $M$ is a real homology pseudosphere, cf. Corollary
\ref{c:folhom}). In the spirit of H. Rummler, \cite{kn:Rum}, and
F. Kamber et al., \cite{kn:KTT}, we obtain lightlike analogs to
Rummler's formula (cf. also \cite{kn:Ton}, p. 66) and to the
transversal divergence theorem, though only on foliated
semi-Riemannian manifolds without boundary, while the problem of
producing a foliated analog of the result by B. \"Unal,
\cite{kn:Una}, is left open.

\vskip 0.1in {\small {\bf Acknowledgements} The present paper was
started while S. Dragomir was a guest of the Department of
Mathematics and Statistics of the University of Windsor (June
2005) and he expresses his gratitude for the excellent working
conditions there. E. Barletta and S. Dragomir were partially
supported by INdAM (Italy) within the interdisciplinary project
{\em Nonlinear subelliptic equations of variational origin in
contact geometry}.}

\section{Screen distributions and transversal bundles}
Let $E \to M$ be a real vector bundle of rank $m$ ($m \geq 2$)
over a $C^\infty$ manifold $M$. In this paper by a ({\em bundle})
{\em metric} in $E$ we intend a $C^\infty$ section $g : x \in M
\mapsto g_x \in E^*_x \otimes_{\mathbb R} E^*_x$ in $E^* \otimes
E^*$ such that $g_x$ is symmetric and has constant index ${\rm
ind}(g_x ) = \sigma$, for any $x \in M$. If each $g_x$ is
nondegenerate and $1 \leq \sigma \leq m-1$ (respectively if each
$g_x$ is positive definite) then $g$ is a {\em semi-Riemannian}
metric (respectively a {\em Riemannian} metric) in $E$. An
arbitrary metric $g$ in $E$ is therefore allowed to be degenerate
i.e. $({\rm Rad} \; E)_x \neq (0)$ for some $x \in M$ where
\[ ({\rm Rad} \; E)_x = \{ v \in E_x : g_x (v , w) = 0, \;\; w \in
E_x \} , \;\;\; x \in M. \] Nevertheless we assume in most cases
that ${\rm Rad} \; E$ is a subbundle of $E$ of rank $r$ with $1
\leq r \leq m$ and then refer to $g$ as a $r$-{\em lightlike}
metric while ${\rm Rad} \; E$ is the {\em radical bundle} of $(E ,
g)$.
\par
Let $\mathcal F$ be a codimension $q$ foliation of a real
$n$-dimensional manifold $M$. Let $\nu ({\mathcal F}) =
T(M)/T({\mathcal F})$ be the transverse bundle and $\Pi : T(M) \to
\nu ({\mathcal F})$ the projection. Let $g$ be a $r$-lightlike
metric in $T({\mathcal F})$ where $1 \leq r \leq \min \{ m , q \}$
and $m = n - q$ . Then $({\mathcal F}, g)$ is a {\em tangentially
lightlike foliation} of $M$ and ${\rm Rad} \; T {\mathcal F}$ is
its {\em tangential radical distribution}. It is with this sort of
foliations that the present paper is mainly concerned. If $M$ is a
$n$-dimensional semi-Riemannian manifold of index $1 \leq s \leq
n-1$ and the metric $g$ above is induced in $T({\mathcal F})$ by
the ambient metric then each leaf of $\mathcal F$ is a degenerate
or lightlike submanifold of $M$ (cf. \cite{kn:DuBe}, p. 140). We
also adopt the terminology in Table 1.
\par
Similarly if $g_Q$ is a $\rho$-lightlike ($1 \leq \rho \leq q$)
metric in $Q = \nu ({\mathcal F})$ such that $\IN{\nabla}_X g_Q =
0$ for any $X \in T({\mathcal F})$ then $({\mathcal F}, g_Q )$ ia
a {\em transversally lightlike foliation}. Here $\IN{\nabla}$
denotes the {\em Bott connection} of $(M , {\mathcal F})$ i.e.
$\IN{\nabla}_X s = \Pi [X, Y]$ for any $C^\infty$ section $s$ in
$Q$ and any $Y \in T(M)$ such that $\Pi (Y) = s$.
\par
Let $\mathcal F$ be a tangentially lightlike foliation of the
semi-Riemannian manifold $(M, g)$. We set
\[ T({\mathcal F})^\bot
 = \{ V \in T(M) : g(V, X) = 0, \;\;\; X \in  T({\mathcal F}) \} . \]
 Let $S(T{\mathcal F})$ and $S(T{\mathcal F}^\bot )$ be complements to the tangential
 radical distribution in $T({\mathcal F})$ and $T({\mathcal
 F})^\bot$, respectively. Then
 \begin{equation}
 T({\mathcal F}) = S(T{\mathcal F}) \oplus {\rm Rad} \; T{\mathcal
 F}, \label{e:fol1}
 \end{equation}
 \begin{equation}
 T({\mathcal F})^\bot = S(T{\mathcal F}^\bot ) \oplus {\rm Rad} \;
 T {\mathcal F} , \label{e:fol2}
 \end{equation}
 and (by Proposition 2.1 in \cite{kn:DuBe}, p. 5) both
 $S(T{\mathcal F})$ and $S(T{\mathcal F}^\bot )$ are
 nondegenerate. Consequently
 \begin{equation}
 T(M) = S(T{\mathcal F}) \oplus S(T{\mathcal F})^\bot .
 \label{e:fol3}
 \end{equation}
 If $T({\mathcal F})$ were nondegenerate then the ambient semi-Riemannian metric
 $g$ would induce a bundle metric
 in $\nu ({\mathcal F})$ by the natural isomorphism $\nu
 ({\mathcal F}) \approx T({\mathcal F})^\bot$. As to the study of
 tangentially lightlike foliations we circumvent the difficulties
 (arising from the failure to decompose $T(M) = T({\mathcal F})
 \oplus T({\mathcal F})^\bot$) by using lightlike transversal
 bundles (as in the theory of lightlike submanifolds, cf. \cite{kn:DuBe}, p.
 139-148).
\[ \begin{array}{|c|c|} \hline \\ {\mathcal F} & r \\ \hline {\rm (I)} \;\;\;
r-{\rm lightlike} & 1 \leq r < \min \{ m, q \}
\\ {\rm (II)} \;\;\; {\rm co-isotropic} & 1 \leq r = q < m \\ {\rm
(III)} \;\;\; {\rm isotropic} & 1 \leq r = m < q \\ {\rm (IV)}
\;\;\; {\rm totally \; lightlike} & 1 \leq r = m = q \\
\hline \end{array} \] \centerline{\small Table 1. Classification
of tangentially lightlike foliations} \centerline{\small of
semi-Riemannian manifolds according to the rank}
\centerline{\small of their tangential radical distribution.}
\vskip 0.1in \noindent Let us start by noticing that
 \begin{equation}
 S(T{\mathcal F})^\bot \supseteq S(T{\mathcal F}^\bot ).
 \label{e:fol4}
 \end{equation}
 Indeed if $X \in S(T{\mathcal F}^\bot ) \subseteq T({\mathcal
 F})^\bot$ then $X$ is orthogonal to $T({\mathcal F}) \supseteq
 S(T{\mathcal F})$ hence $X \in S(T{\mathcal F})^\bot$. Next,
 since $S(T{\mathcal F}^\bot )$ is nondegenerate
 \begin{equation}
 S(T{\mathcal F})^\bot = S(T{\mathcal F}^\bot ) \oplus
 S(T{\mathcal F}^\bot )^\bot . \label{e:fol5}
 \end{equation}
 We shall need the following
 \begin{lemma} Let $\{ \xi_1 , \cdots , \xi_r \}$ be a local frame
 of ${\rm Rad} \; T{\mathcal F}$ defined on the open set $U
 \subseteq M$. There exist $N_i \in \Gamma^\infty (U ,
 S(T{\mathcal F}^\bot )^\bot )$, $1 \leq i \leq r$, such that
 $g(N_i , \xi_j ) = \delta_{ij}$ and $g(N_i , N_j ) = 0$.
 \label{l:fol1}
 \end{lemma}
{\em Proof}. Note first that
\begin{equation}
{\rm Rad} \; T{\mathcal F} \subseteq S(T{\mathcal F}^\bot )^\bot .
\label{e:fol6}
\end{equation}
Indeed if $X \in {\rm Rad} \; T{\mathcal F}$ then $X$ is
orthogonal on $T({\mathcal F})^\bot \supseteq S(T{\mathcal F}^\bot
)$ hence $X \in S(T{\mathcal F}^\bot )^\bot$. Next we choose a
complement $E$ to ${\rm Rad} \; T{\mathcal F}$ so that
\begin{equation}
S(T{\mathcal F}^\bot )^\bot = ({\rm Rad} \; T{\mathcal F}) \oplus
E. \label{e:fol7}
\end{equation}
Consequently $\dim_{\mathbb R} E_x = r$ for any $x \in M$. Let
then $\{ V_1 , \cdots , V_r \}$ be a local frame of $E$ on $U$.
One may look for the $N_i$'s in the form
\[ N_i = A_i^k \xi_k + B_i^k V_k \]
for some $C^\infty$ functions $A_i^k , B_i^k : U \to {\mathbb R}$
with the requirement
\[ \delta_{ij} = g(N_i , \xi_j ) = B_i^k g_{jk} \]
where $g_{jk} = g(\xi_j , V_k )$. Let us set $G = \det [g_{jk}]$.
We claim that $G(x) \neq 0$ for any $x \in U$. The proof is by
contradiction. If $G(x_0 ) = 0$ for some $x_0 \in U$ then there is
$v = (v^1 , \cdots , v^r ) \in {\mathbb R}^r \setminus \{ 0 \}$
such that
\begin{equation}
g_{jk}(x_0 ) v^j = 0, \;\;\; 1 \leq k \leq r.
\label{e:fol8}
\end{equation}
Let us set $w = v^j \xi_{j,x_0} \in ({\rm Rad} \; T{\mathcal
F})_{x_0} \subset S(T{\mathcal F}^\bot )^\bot$. Then (by
(\ref{e:fol8})) $g_{x_0} (w , V_{k,x_0} ) = 0$. Also $g_{x_0} (w ,
\xi_{k,x_0}) = 0$ by the very definition of $w$. Then $w$ sits in
$S(T{\mathcal F}^\bot )^\bot_{x_0} \setminus \{ 0 \}$ and (by
(\ref{e:fol7})) it is perpendicular on $S(T{\mathcal F}^\bot
)^\bot_{x_0}$ i.e. $S(T{\mathcal F}^\bot )^\bot_{x_0}$ is
degenerate, a contradiction. Therefore it is legitimate to
consider $[g^{jk}] := [g_{jk}]^{-1}$. Then $B^k_i = g^{ki}$ and
the requirement $g(N_i , N_j ) = 0$ yields
\[ A^j_i + A_j^i + g^{ki} g^{\ell j} g(V_k , V_\ell ) = 0 \]
and we may choose $A^i_j := - \frac{1}{2} g^{ki} g^{\ell j} g(V_k
, V_\ell )$. Lemma \ref{l:fol1} is proved. In particular (with the
notations of Lemma \ref{l:fol1}) $\{ \xi_1 , \cdots , \xi_r , N_1
, \cdots , N_r \}$ is a local frame of $S(T{\mathcal F}^\bot
)^\bot$ on $U$. Let us set
\[ ltr(T{\mathcal F})_x = \sum_{i=1}^r {\mathbb R} N_{i,x} \, ,
\;\;\; x \in U. \]
\begin{lemma} The definition of the bundle $ltr(T{\mathcal F})_x$ doesn't
depend upon the choice of local frames $\{ \xi_j \}$ of ${\rm Rad}
\; T{\mathcal F}$ and $\{ V_k \}$ of $E$ at $x$. Moreover
$ltr(T{\mathcal F}) = \cup_{x \in M} ltr (T{\mathcal F})_x$ is a
vector bundle over $M$ and
\begin{equation}
S(T{\mathcal F}^\bot )^\bot = ({\rm Rad} \; T {\mathcal F}) \oplus
ltr (T{\mathcal F}). \label{e:fol9}
\end{equation}
\label{l:folomit}
\end{lemma}
The proof of Lemma \ref{l:folomit} is imitative of that of Theorem
1.4 in \cite{kn:DuBe}, p. 147, and is omitted. We call
$ltr(T{\mathcal F}) \to M$ a {\em lightlike transversal} vector
bundle with respect to the pair $(S(T{\mathcal F}), \;
S(T{\mathcal F}^\bot ))$. Also
\begin{equation}
tr(T{\mathcal F}) := ltr(T{\mathcal F}) \oplus S(T{\mathcal
F}^\bot ) \label{e:fol10}
\end{equation}
is a {\em transversal} vector bundle. Then (by (\ref{e:fol3}),
(\ref{e:fol5}) and (\ref{e:fol9}))
\[ T(M) = S(T{\mathcal F}) \oplus S(T{\mathcal F}^\bot ) \oplus
({\rm Rad}\; T{\mathcal F}) \oplus ltr(T{\mathcal F}) \] hence
\begin{equation}
T(M) = T({\mathcal F}) \oplus tr(T{\mathcal F}). \label{e:fol11}
\end{equation}
Let $\sigma : \nu ({\mathcal F}) \to tr(T{\mathcal F})$ be the
bundle isomorphism given by
\[ \sigma (s) = {\rm tra} (Y), \;\;\; \Pi (Y) = s, \;\;\; Y \in
T(M), \] where ${\rm tra} : T(M) \to tr(T{\mathcal F})$ is the
natural projection associated to the decomposition
(\ref{e:fol11}). Let us set
\[ g_{\rm tra} (s, r) = g(\sigma (s) \, , \, \sigma (r)), \;\;\;
s,r \in \nu ({\mathcal F}). \] If $g_{\rm tra}$ is holonomy
invariant, i.e. ${\mathcal L}_X g_{\rm tra} = 0$ for any $X \in
T({\mathcal F})$, then $g$ is said to be {\em bundle-like}. Here
${\mathcal L}_X$ denotes the Lie derivative in the direction $X$.
Let $Q = \nu ({\mathcal F})$ for simplicity. One expects $g_{\rm
tra}$ to be degenerate, as well. Indeed, if we set
\[ {\rm Rad} \; Q = \{ s \in Q : g_{\rm tra}(s,r) = 0, \;\; r \in
Q \}  \] then we have
\begin{proposition} Let $\mathcal F$ be a lightlike foliation of the
semi-Riemannian manifold $(M , g)$ and $ltr(T{\mathcal F}) \to M$
a lightlike transversal vector bundle associated with the screen
distributions $S(T{\mathcal F})$ and $S(T{\mathcal F}^\bot )$.
Then
\begin{equation}
\sigma ({\rm Rad} \; Q) = ltr(T{\mathcal F}). \label{e:fol12}
\end{equation} \label{p:fol1}
\end{proposition}
{\em Proof}. Let $N \in ltr(T{\mathcal F})$ and $r \in Q$. As $N$
is orthogonal to $tr(T{\mathcal F})$
\[ g_{\rm tra}(\sigma^{-1} (N) , r) = g(N , \sigma (r)) = 0 \]
it follows that $\sigma^{-1} (N) \in  {\rm Rad} \; Q$. For the
opposite inclusion let $s \in {\rm Rad} \; Q$ and $Z \in
tr(T{\mathcal F})$. If we set $r = \sigma^{-1} (Z) \in Q$ then
\begin{equation} 0 = g_{\rm tra} (s,r) = g(\sigma (s) , Z).
\label{e:fol13}
\end{equation}
We have $s = \Pi (Y)$ for some $Y \in T(M)$. As a consequence of
(\ref{e:fol10})-(\ref{e:fol11}) $Y = X + N + V$ for some $X \in
T({\mathcal F})$, $N \in ltr(T{\mathcal F})$ and $V \in
S(T{\mathcal F}^\bot )$. Then $\sigma (s) = N + V$. Let $W \in
S(T{\mathcal F}^\bot )$. Applying (\ref{e:fol13}) for $Z = W$
gives $g(V, W) = 0$ and then $V = 0$ since $S(T{\mathcal F}^\bot
)$ is nondegenerate. It remains that $\sigma (s) = N \in {\rm Rad}
\; Q$ and Proposition \ref{p:fol1} is proved.
\par
By (\ref{e:fol10}) and Proposition \ref{p:fol1} a canonical choice
of screen distribution in $Q$ is $S(Q) := \sigma^{-1} S(T{\mathcal
F}^\bot )$ so that
\[ Q = S(Q) \oplus {\rm Rad} \; Q . \]
Let $\nabla^g$ be the Levi-Civita connection of $(M , g)$. We also
consider
\[ \nabla_X s = \begin{cases} \IN{\nabla}_X s , & X \in
T({\mathcal F}), \cr \Pi \nabla^g_X \sigma (s), & X \in
ltr(T{\mathcal F}). \cr \end{cases} \] One checks easily that
\begin{proposition}
$\nabla$ is a connection in $Q$ and $T_\nabla = 0$, where
$T_\nabla (Y,Z) := \nabla_Y \Pi Z - \nabla_Z \Pi Y - \Pi [Y,Z]$
for any $Y,Z \in T(M)$. Moreover $g$ is bundle-like if and only if
$\nabla g_{\rm tra} = 0$.
\end{proposition}

\section{Lightlike foliations defined by suspension} Let $(N, h)$
be a semi-Riemannian manifold and $j : B \hookrightarrow N$ a
$m$-dimensional connected lightlike submanifold i.e. $({\rm Rad}
\; TB)_x = \{ X \in T_x (B) : g_{B,x}(X,Y) = 0, \;\; Y \in T_x (B)
\}$ $(x \in B)$ has constant dimension $1 \leq \rho \leq \min \{
m, \; \ell \}$ ($\ell = \dim_{\mathbb R} N - m$). Here $g_B = j^*
h$. Let $\tilde{B}$ be the universal covering manifold of $B$ and
$\tilde{p} : \tilde{B} \to B$ the projection. We set
$g_{\tilde{B}} = \tilde{p}^* g_B$. Next let us consider a
$q$-dimensional connected semi-Riemannian manifold $(T, g_T )$ and
the {\em warped product} $\hat{M} = T \times_f \tilde{B}$ i.e.
$\hat{M}$ is the product manifold $T \times \tilde{B}$ endowed
with the $(0,2)$-tensor field $\hat{g} = p_1^* \; g_T + (f \circ
p_1 )^2 \; p_2^* \; g_{\tilde{B}}$ where $f$ (the {\em warping
function}) is a $C^\infty$ function $f : T \to (0 , + \infty )$.
Also $p_1 : \hat{M} \to T$ and $p_2 : \hat{M} \to \tilde{B}$ are
the natural projections. Our notion of warped product generalizes
slightly that in \cite{kn:O'Nei}, p. 204, as $\hat{g}$ is not a
semi-Riemannian metric ($(\hat{M}, \hat{g})$ has a nontrivial
radical distribution). Let $\hat{\mathcal F}$ be the foliation of
$\hat{M}$ whose leaves are the fibres of $p_1$ i.e.
$T(\hat{\mathcal F}) = {\rm Ker}(d p_1 )$ and
$\hat{M}/\hat{\mathcal F} = \{ \{ y \} \times \tilde{B} : y \in T
\}$. The same symbol $\hat{g}$ denotes the induced metric in
$T(\hat{\mathcal F})$.
\begin{lemma} $\hat{g}$ is a $\rho$-lightlike metric in
$T(\hat{\mathcal F})$.
\end{lemma}
\noindent {\em Proof}. We set as customary
\[ ( {\rm Rad} \; T \hat{\mathcal F} )_{(y, \tilde{x})}
= \{ X \in T(\hat{\mathcal F})_{(y , \tilde{x})} : \hat{g}_{(y,
\tilde{x})} (X,Y) = 0, \;\; Y \in T(\hat{\mathcal F})_{(y,
\tilde{x})} \} \] for any $y \in T$ and $\tilde{x} \in \tilde{B}$.
If $\alpha_y : \tilde{B} \to \hat{M}$ is the canonical injection
$\alpha_y (\tilde{x}) = (y , \tilde{x})$ then
\begin{equation}
({\rm Rad} \; T\hat{\mathcal F})_{(y, \tilde{x})} = (d_{\tilde{x}}
\alpha_y ) ({\rm Rad} \; T \tilde{B})_{\tilde{x}} \, .
\label{e:fol15}
\end{equation}
Indeed
\[ ({\rm Rad} \; T\hat{\mathcal F})_{(y, \tilde{x})} = \{ X \in
{\rm Ker}(d_{(y, \tilde{x})} p_1 ) : \left( p_1^* \; g_T + \right.
\]
\[ \left. + (f \circ p_1 )^2 \;\; p_2^* \; g_{\tilde{B}} \right)_{(y ,
\tilde{x})} (X,Y) = 0, \;\; Y \in {\rm Ker}(d_{(y, \tilde{x})} p_1
) \} = \]
\[ = \{ X \in
{\rm Ker}(d_{(y, \tilde{x})} p_1 ) : g_{\tilde{B} , \tilde{x}}
((d_{(y, \tilde{x})} p_2 ) X , (d_{(y, \tilde{x})} p_2 ) Y ) = 0,
\] \[  Y \in {\rm Ker}(d_{(y, \tilde{x})} p_1 ) \} .
 \]
Note that ${\rm Ker}(d_{(y, \tilde{x})} p_1 ) = (d_{\tilde{x}}
\alpha_y ) T_{\tilde{x}} (\tilde{B})$. Let $X = (d_{\tilde{x}}
\alpha_y ) v$ with $v \in T_{\tilde{x}} (\tilde{B})$. As $p_2
\circ \alpha_y = 1$ (the identical transformation of $\tilde{B}$)
it follows that $(d_{(y , \tilde{x})} p_2 ) X = v$. We conclude
that
\[ ({\rm Rad} \; T\hat{\mathcal F})_{(y, \tilde{x})} = \{
(d_{\tilde{x}} \alpha_y ) v : g_{\tilde{B} , \tilde{x}} (v , w) =
0, \;\; w \in T_{\tilde{x}} (\tilde{B}) \} \] and (\ref{e:fol15})
is proved. $\square$
\par
Let $\varphi : \pi_1 (B, x_0 ) \to {\rm Diff}(T)$ be a
homomorphism of the fundamental group of $B$ with base point $x_0
\in B$ into the group of all $C^\infty$ diffeomorphisms of $T$ in
itself. Also we think of $\tilde{B}$ as the set of all homotopy
classes of paths issuing at $x_0$. Let $\hat{M} \times \pi_1 (B ,
x_0 ) \to \hat{M}$ be the natural action given by
\[ R_{[\gamma ]} (y , \tilde{x}) = (\varphi ([\gamma ]^{-1} )(y)
\, , \, \tilde{x} \cdot [\gamma ]), \;\;\; (y, \tilde{x}) \in
\hat{M}, \;\; [\gamma ] \in \pi_1 (B , x_0 ), \] and let $M =
\hat{M} /\pi_1 (B , x_0 )$ be the quotient space. Let $\mathcal F$
be the projection of $\hat{\mathcal F}$ on $M$ i.e. the foliation
of $M$ whose leaves are the projection of the leaves of
$\hat{\mathcal F}$
\[ M/{\mathcal F} = \{ \pi (L) : L \in \hat{M}/\hat{\mathcal F} \}
\]
where $\pi : \hat{M} \to M$ denotes the natural projection.
Infinitesimally $T({\mathcal F})_{\pi (y, \tilde{x})} = (d_{(y,
\tilde{x})} \pi ) T(\hat{\mathcal F})_{(y, \tilde{x})}$. We say
$\mathcal F$ is the foliation of $M$ defined by {\em suspension}
of the homomorphism $\varphi$. Let us consider the map $p : M \to
B$ given by $p(\pi (y , \tilde{x})) = \tilde{p}(\tilde{x})$. Then
$p : M \to B$ is a fibre bundle with standard fibre $T$ and
structure group $G = \varphi (\pi_1 (B , x_0 ))$. See also
\cite{kn:Mol}, p. 29. As well known (cf. e.g. \cite{kn:Mol}, p.
28) $M$ admits a natural $C^\infty$ manifold structure such that
$\pi : \hat{M} \to M$ is an \'etale mapping i.e. ${\rm Ker}(d_{(y,
\tilde{x})} \pi ) = (0)$ for any $(y , \tilde{x}) \in \hat{M}$.
\begin{lemma} If $G \subset {\rm Isom}(T)$ and $f$ is
$G$-invariant then $\hat{g}$ is $\pi_1 (B , x_0 )$-invariant. In
particular there is a metric $g$ in $T(M)$ such that $\pi^* g =
\hat{g}$.
\label{l:folex4}
\end{lemma}
\noindent Here ${\rm Isom}(T)$ denotes the group of isometries of
the semi-Riemannian manifold $(T , g_T )$. {\em Proof of Lemma}
\ref{l:folex4}. Let $[\gamma ] \in \pi_1 (B , x_0 )$ and $u, v \in
T_{(y, \tilde{x})} (\hat{M})$. Then
\[ \hat{g}_{R_{[\gamma ]}(y , \tilde{x})} (A , B)
= g_{T, \varphi ([\gamma ]^{-1})(y)} ((d_{R_{[\gamma ]}(y,
\tilde{x})} p_1 ) A \, , \, (d_{R_{[\gamma ]}(y, \tilde{x})} p_1 )
B ) + \]
\[ + f(\varphi ([\gamma ]^{-1})(y))^2 g_{\tilde{B}, \tilde{x}
\cdot [\gamma ]} ((d_{R_{[\gamma ]}(y , \tilde{x})} p_2 ) A \, ,
\, (d_{R_{[\gamma ]}(y , \tilde{x})} p_2 ) B) \] where $A =
(d_{(y, \tilde{x})} R_{[\gamma ]}) u$ and $B = (d_{(y, \tilde{x})}
R_{[\gamma ]}) v$. Note that
\[ p_1 \circ R_{[\gamma ]} = \varphi ([\gamma ]^{-1}) \circ p_1 \,
, \;\;\; p_2 \circ R_{[\gamma ]} = D_{[\gamma ]} \circ p_2 \, , \]
where $D_{[\gamma ]} : \tilde{B} \to \tilde{B}$ is the deck
transformation $D_{[\gamma ]}(\tilde{x}) = \tilde{x} \cdot [\gamma
]$. Then
\[ (R^*_{[\gamma ]} \hat{g})_{(y, \tilde{x})} (u,v) = (\varphi
([\gamma ]^{-1})^* g_T )_y ((d_{(y, \tilde{x})} p_1 ) u \, , \,
(d_{(y, \tilde{x})} p_1 ) v) + \]
\[ + f(\varphi ([\gamma ]^{-1})(y))^2 (D^*_{[\gamma ]}
g_{\tilde{B}} )_{\tilde{x}} ((d_{(y , \tilde{x})} p_2 ) u \, , \,
(d_{(y, \tilde{x})} p_2 ) v) = \hat{g}_{(y, \tilde{x})} (u,v). \]
To prove the second statement in Lemma \ref{l:folex4} let $p \in
M$ and $X,Y \in T_p (M)$. Then $p = \pi (y, \tilde{x})$ and $X =
(d_{(y, \tilde{x})} \pi ) u$, $Y = (d_{(y , \tilde{x})} \pi ) v$
for some $y \in T$, $\tilde{x} \in \tilde{B}$ and $u, v \in T_{(y,
\tilde{x})} (\hat{M})$. We set
\[ g_p (X,Y) := \hat{g}_{(y, \tilde{x})} (u , v). \]
We only need to check that the definition doesn't depend upon the
choice of representatives. If $(y^\prime , \tilde{x}^\prime ) \in
\pi^{-1} (p)$ then $(y^\prime , \tilde{x}^\prime ) = R_{[\gamma
]}(y , \tilde{x})$ for some $[\gamma ] \in \pi_1 (B , x_0 )$. If
$u^\prime , v^\prime \in T_{(y^\prime , \tilde{x}^\prime )}
(\hat{M})$ are other representatives of $X$ and $Y$ then $\pi =
\pi \circ R_{[\gamma ]}$ yields $u^\prime - (d_{(y, \tilde{x})}
R_{[\gamma ]})u \in {\rm Ker}(d_{(y , \tilde{x})} \pi ) = (0)$
i.e. $u^\prime = (d_{(Y, \tilde{x})} R_{[\gamma ]} ) u$ and
similarly $v^\prime = (d_{(y, \tilde{x})} R_{[\gamma ]} ) v$.
Finally (by Lemma \ref{l:folex4}) $\hat{g}_{(y^\prime ,
\tilde{x}^\prime )} (u^\prime , v^\prime ) = (R^*_{[\gamma ]}
\hat{g})_{(y , \tilde{x})} (u, v) = \hat{g}_{(y, \tilde{x})}
(u,v)$ i.e. $g_p (X,Y)$ is well defined. $\square$
\par
The same symbol $g$ denotes the induced metric in $T({\mathcal
F})$.
\begin{proposition} Let $T$ be a connected semi-Riemannian manifold, $j : B
\hookrightarrow N$ a connected lightlike submanifold of the
semi-Riemannian manifold $(N, h)$ such that the induced metric
$j^* h$ is $\rho$-lightlike, and $\tilde{B}$ the universal
covering manifold of $B$. Let $\varphi : \pi_1 (B , x_0 ) \to {\rm
Isom}(T)$ be a group homomorphism and let $\hat{M} = T \times_f
\tilde{B}$ be a warped product with a $G$-invariant warping
function $f : T \to (0, + \infty )$ where $G = \varphi (\pi_1 (B ,
x_0 ))$. Let $M = \hat{M}/\pi_1 (B , x_ 0)$. Then the foliation
$\mathcal F$ of $M$ defined by suspension of the homomorphism
$\varphi$ is tangentially lightlike and $g$ is a $\rho$-lightlike
metric in $T({\mathcal F})$.
\end{proposition}
{\em Proof}. The tangential radical distribution is $({\rm Rad} \;
T {\mathcal F} )_{\pi (y , \tilde{x})} =$ \[ = \{ X \in
T({\mathcal F})_{\pi (y, \tilde{x})} : g_{\pi (y, \tilde{x})}
(X,Y) = 0, \;\; Y \in T({\mathcal F})_{\pi (y, \tilde{x})} \} =
\]
\[ = \{ (d_{(y, \tilde{x})} \pi ) u : (\pi^* g)_{(y, \tilde{x})}
(u, v) = 0, \;\; u,v \in T(\hat{\mathcal F})_{(y, \tilde{x})} \} =
\] \[ = (d_{(y, \tilde{x})} \pi ) ({\rm Rad} \; T \hat{\mathcal
F})_{(y, \tilde{x})} = (d_{(y, \tilde{x})} \pi ) (d_{\tilde{x}}
\alpha_y ) ({\rm Rad} \; T \tilde{B} )_{\tilde{x}}  \] (by
(\ref{e:fol15})). Next $\tilde{p}^* g_B = g_{\tilde{B}}$ yields
\begin{equation}
(d_{\tilde{x}} \tilde{p}) ({\rm Rad} \; T \tilde{B})_{\tilde{x}} =
({\rm Rad} \; T B)_{\tilde{p}(\tilde{x})} \label{e:folex15}
\end{equation}
for any $\tilde{x} \in \tilde{B}$. Our previous calculation, the
identity (\ref{e:folex15}) and the commutativity of the diagram
\[ \begin{array}{ccccc} \empty & \tilde{B} &
\stackrel{\alpha_y}{\longrightarrow} & \hat{M} & \empty \\
\tilde{p} & \downarrow & \empty & \downarrow & \pi \\ \empty & B &
\stackrel{p}{\longleftarrow} & M & \empty \end{array} \] imply
that
\begin{equation}
(d_{\pi (y , \tilde{x})} p) ({\rm Rad} \; T {\mathcal F})_{\pi (y
, \tilde{x})} = ({\rm Rad}\; T B)_{\tilde{p}(\tilde{x})} \, .
\end{equation}
Cf. again \cite{kn:Mol}, p. 29, the fibres of $p$ are connected
total transversals of $(M , {\mathcal F})$. In particular
\[ T_{\pi (y , \tilde{x})} (M) = T({\mathcal F})_{\pi (y,
\tilde{x})} \oplus {\rm Ker}(d_{\pi (y , \tilde{x})} p), \;\;\;
\pi (y, \tilde{x}) \in p^{-1} (\tilde{p}(\tilde{x})). \] Hence the
restriction of $d_{\pi (y , \tilde{x})} p$ to $({\rm Rad} \; T
{\mathcal F})_{\pi (y , \tilde{x})}$ is a ${\mathbb R}$-linear
isomorphism.

\section{Lightlike Killing vector fields}
For each lightlike foliation $\mathcal F$ of the semi-Riemannian
manifold $M$ we denote by $\Omega^k_B ({\mathcal F})$ the space of
all basic differential $k$-forms on $(M, {\mathcal F})$ i.e. if
$\omega \in \Omega^k_B ({\mathcal F})$ then $X \, \rfloor \,
\omega = 0$ and $X \, \rfloor \, d \omega = 0$ for any $X \in
T({\mathcal F})$. In particular $\Omega^0_B ({\mathcal F})$ is the
space of all basic functions ($f \in C^\infty (M)$ is basic if $X
(f) = 0$ for any $X \in T({\mathcal F})$). Let
\[ H^k_B ({\mathcal F}) := H^k (\Omega^\bullet_B ({\mathcal F})),
\;\;\; k \geq 0, \] be the corresponding cohomology groups (that
is the {\em basic cohomology} of $(M , {\mathcal F}))$. By
standard foliation theory $H^0_B ({\mathcal F}) = {\mathbb R}$ and
there is a natural injection $H^1_B ({\mathcal F}) \hookrightarrow
H^1 (M , {\mathbb R})$ (cf. e.g. \cite{kn:Ton}, p. 119).
\par
Let $(M , g)$ be a geodesically complete semi-Riemannian manifold
and $\xi$ a lightlike Killing vector field on $M$. Then (cf. e.g.
\cite{kn:O'Nei}, p. 254) $\xi$ is complete. Let $H$ be the global
$1$-parameter group of global transformations of $M$ obtained by
integrating $\xi$. Let $G = \overline{H}$ be the closure of $H$ in
${\rm Isom}(M , g)$. We assume from now on that $G$ is compact.
For instance, if $(M , g)$ is a Lorentz manifold ($s = 1$) and
${\rm Isom}(M, g) = {\rm O}(n+1,1)$ then the closure of any
lightlike $1$-parameter subgroup is compact, cf. Lemma 3.1 in
\cite{kn:Kam}, p. 584. Let $\Omega^\bullet (M)$ be the de Rham
algebra of $M$ and $\Omega^\bullet (M)^G$ the subalgebra of all
$G$-invariant differential forms i.e. if $\omega \in
\Omega^\bullet (M)^G$ then ${\mathcal L}_\xi \omega = 0$. Let
${\mathcal F}$ be the codimension $q = n-1$ lightlike foliation of
$M$ such that $T({\mathcal F}) = {\mathbb R} \xi$. It is immediate
that
\begin{proposition}
Either $\mathcal F$ is isotropic or $M$ is a Lorentz surface {\rm
(}$n = 2$, $s = 1${\rm )} and $\mathcal F$ is totally lightlike.
\end{proposition}
Next, note that
\begin{equation} \Omega^k_B ({\mathcal F}) \subset \Omega^k (M)^G
\, , \;\;\; k \geq 0.
\label{e:fol14}
\end{equation}
Let $i_\xi : \Omega^k (M) \to \Omega^{k-1}(M)$ be the interior
product with $\xi$ i.e. $i_\xi \omega = \xi \, \rfloor \, \omega$
for any $\omega \in \Omega^k (M)$. Then
\begin{equation}
i_\xi \; \Omega^k (M)^G \subseteq \Omega_B^{k-1}({\mathcal F}),
\;\;\; k \geq 1. \label{:fol15}
\end{equation}
Indeed, let $\omega \in \Omega^k (M)^G$ and $\eta = i_\xi \omega$.
Then $i_\xi \eta = i_\xi^2 \omega = 0$ and
\[ i_\xi d \eta = i_\xi {\mathcal L}_\xi \omega - i_\xi^2 d \omega
= 0 \] by Cartan's formula. Our main purpose in the present
section is to establish
\begin{theorem} $\,$ \par\noindent Let $\xi$ be a lightlike Killing vector field on the
complete semi-Riemannian manifold $(M,g)$ and $\mathcal F$ the
$1$-dimensional foliation tangent to $\xi$. Let $G$ be the closure
in ${\rm Isom}(M, g)$ of the $1$-parameter group generated by
$\xi$. Assume that there is a globally defined $G$-invariant
vector field $V \neq 0$ on $M$ such that
\begin{equation}
\label{e:fol16} S(T{\mathcal F}^\bot )^\bot = ({\rm Rad} \; T
{\mathcal F}) \oplus {\mathbb  R} V.
\end{equation}
If $G$ is compact then for any $k \geq 1$ there is a linear map
$\Delta : H^{k-1}_B ({\mathcal F}) \to H^{k+1}_B ({\mathcal F})$
such that
\begin{equation}
\label{e:fol17} H^k_B ({\mathcal F})
\stackrel{j_*}{\longrightarrow} H^k (M , {\mathbb R})
\stackrel{(i_\xi )_*}{\longrightarrow} H^{k-1}_B ({\mathcal F})
\stackrel{\Delta}{\longrightarrow} H^{k+1}_B ({\mathcal F}) \to
\cdots
\end{equation}
is a long exact sequence, where $j : \Omega^k_B ({\mathcal F}) \to
\Omega^k (M)^G$ is the inclusion. In particular, if $M$ is compact
then $\dim_{\mathbb R} H^k_B ({\mathcal F}) < \infty$.
\label{t:fol1}
\end{theorem}
{\em Proof}. The map $i_\xi : \Omega^k (M)^G \to \Omega_B^{k-1}
({\mathcal F})$ is surjective. Indeed, let $V \in S(T{\mathcal
F}^\bot )$ as in Theorem \ref{t:fol1}. Then (by the proof of Lemma
\ref{l:fol1}) $g(\xi , V) \neq 0$ everywhere on $M$. Let us set
\[ N = \frac{1}{g(\xi , V)} \{ V - \frac{g(V,V)}{2 g(\xi , V)} \;
\xi \}  \] so that $g(\xi , N) = 1$ and $g(N,N) = 0$. Since
\[ \xi (g(\xi , V)) = ({\mathcal L}_\xi g)(\xi , V) + g(\xi ,
{\mathcal L}_\xi V) = 0, \]
\[ \xi (g(V,V)) = ({\mathcal L}_\xi g)(V,V) + 2 g({\mathcal L}_\xi
V , V) = 0, \] it follows that
\begin{equation}
{\mathcal L}_\xi N = 0. \label{e:fol18}
\end{equation}
Let us consider the $1$-form $\alpha \in \Omega^1 (M)$ given by
$\alpha (X) = g(X,N)$ for any $X \in T(M)$. Note that $\alpha$ is
$G$-invariant. Indeed (by (\ref{e:fol18}))
\[ ({\mathcal L}_\xi \alpha ) X = \xi (g(X,N)) - g({\mathcal
L}_\xi X , N) = ({\mathcal L}_\xi g)(X, N) = 0 \] for any $X \in
T(M)$. Consequently, for any $\omega \in \Omega^{k-1}_B ({\mathcal
F})$ the $k$-form $\alpha \wedge \omega$ is $G$-invariant. Finally
\[ i_\xi (\alpha \wedge \omega ) = (i_\xi \alpha ) \omega = \omega
\]
so that $i_\xi$ is on-to, as claimed. The next step is to observe
that
\begin{equation}
\label{e:fol19} 0 \to \Omega_B^k ({\mathcal F}) \to \Omega^k (M)^G
\stackrel{i_\xi}{\longrightarrow} \Omega^{k-1}_B ({\mathcal F})
\to 0
\end{equation}
is a short exact sequence. By (\ref{e:fol14}) and the first part
of the proof of Theorem \ref{t:fol1} one only needs to check
exactness at the middle term. If $\omega \in \Omega^k_B ({\mathcal
F})$ then $i_\xi \omega = 0$ because $\omega$ is a basic form.
Viceversa, let $\omega \in {\rm Ker}(i_\xi ) \subseteq \Omega^k
(M)^G$. Then $i_\xi \omega = 0$ and ${\mathcal L}_\xi \omega = 0$
hence $\omega \in \Omega^k_B ({\mathcal F})$.
\par
Let us consider the map \[ \Delta : H^{k-1}_B ({\mathcal F}) \to
H^{k+1}_B ({\mathcal F}) \] given by $\Delta [\omega ] = [d \alpha
\; \wedge \; \omega ]$ for any $\omega \in \Omega^{k-1}_B
({\mathcal F})$ with $d \omega = 0$. As $i_\xi \omega = 0$ and
$i_\xi \alpha = 1$ one has
\[ i_\xi (d \alpha \; \wedge \; \omega ) = (i_\xi \; d \alpha
)\wedge \omega = ({\mathcal L}_\xi \alpha - d \; i_\xi \alpha )
\wedge \omega = 0, \]
\[ {\mathcal L}_\xi (d \alpha \; \wedge \; \omega ) = ({\mathcal
L}_\xi \; d \alpha )\wedge \omega = (d \; i_\xi \; d \alpha )
\wedge \omega = 0, \] hence $d \alpha \; \wedge \; \omega \in
\Omega^{k+1}_B ({\mathcal F})$. Also the form $d \alpha \; \wedge
\; \omega$ is closed, so that its cohomology class $\bmod \; d \;
\Omega_B^k ({\mathcal F})$ is well defined. One checks easily that
the definition of $\Delta [\omega ]$ doesn't depend upon the
choice of representative in $[\omega ]$. At this point one may use
the sequence (\ref{e:fol19}) and the map $\Delta$ to build the
sequence
\[ H^k_B ({\mathcal F}) \to H^k (\Omega^\bullet (M)^G ) \to
H^{k-1}_B ({\mathcal F}) \to H^{k+1}_B ({\mathcal F}) \to \cdots
\]
which yields (\ref{e:fol17}) as the compactness of $G$ implies
$H^k (\Omega^\bullet (M)^G ) \approx H^k (M , {\mathbb R})$ (cf.
e.g. \cite{kn:GHV}, p. 151). Since (\ref{e:fol19}) is already
exact we need to check exactness in (\ref{e:fol17}) only at the
terms of the form $H^{k-1}_B ({\mathcal F})$. For any $\omega \in
\Omega^k (M)^G$ with $d \omega = 0$ we have
\[ \Delta (i_\xi )_* [\omega ] = [d \alpha \; \wedge \; i_\xi
\omega ] = [d(\alpha \wedge i_\xi \omega ) + \alpha \wedge d i_\xi
\omega ] = \]
\[ = [\alpha \wedge ({\mathcal L}_\xi \omega - i_\xi d \omega )] =
0. \] Viceversa, if $[\eta ] \in {\rm Ker}(\Delta )$ then $\eta
\in \Omega_B^{k-1} ({\mathcal F})$ and $d \eta = 0$ and
\[ d \alpha \; \wedge \; \eta =  d \beta \]
for some $\beta \in \Omega^k_B ({\mathcal F})$. Then $\omega :=
\alpha \wedge \eta -\beta$ is a closed $G$-invariant form and
$i_\xi \beta = 0$ yields $(i_\xi )_* [\omega ] = [\eta ]$. Theorem
\ref{t:fol1} is proved. $\square$ \vskip 0.1in With the notations
above a lightlike Killing vector field $\xi$ is said to be {\em
complemented} if there exist nowhere zero globally defined
$G$-invariant vector fields $W \in S(T{\mathcal F}^\bot )$ and $N
\in ltr(T{\mathcal F})$ such that 1) $[W,N] = 0$, 2) $W$ is
spacelike and $g(\xi , N) = 1$, and 3) for any $f \in \Omega^0_B
({\mathcal F})$ there are $a,b \in \Omega^0_B ({\mathcal F})$ such
that $N(b) - W(a) = f$.
\begin{proposition} Let $M$ be a $3$-dimensional semi-Riemannian
manifold and $\xi$ a complemented Killing vector field on $M$.
Then $H^2_B ({\mathcal F}) = 0$. \label{p:fol3}
\end{proposition}
{\em Proof}. We may assume without loss of generality that $g(W ,
W) = 1$. Otherwise we set $W^\prime = g(W,W)^{-1/2} W$ and observe
that ${\mathcal L}_\xi g = 0$ and ${\mathcal L}_\xi W = 0$ yield
${\mathcal L}_\xi W^\prime = 0$. Let us set
\[ \lambda (X) = g(X,N), \;\; \mu (X) = g(X, \xi ), \;\; \eta (X)
= g(X,W), \] for any $X \in T(M)$. Then any $\Omega \in \Omega^2_B
({\mathcal F})$ is given by $\Omega = f \; \mu \wedge \eta$ for
some $f \in \Omega_B^0 ({\mathcal F})$. Indeed $\xi \, \rfloor \,
\Omega = 0$ implies that $\Omega^2_B ({\mathcal F})$ is spanned by
$\mu \wedge \eta$. To see that the coefficient is a basic function
one must compute $\xi \, \rfloor \, d \Omega$. Note that (by $\mu
(\xi ) = 0$)
\[ 2 (\xi \, \rfloor \, d \mu ) X = \xi (\mu (X)) -
X(\mu (\xi )) - \mu ([\xi , X]) = \] \[ = \xi (g(X, \xi )) -
g({\mathcal L}_\xi X , \xi ) = ({\mathcal L}_\xi g)(X, \xi ) = 0
\] for any $X \in T(M)$. Therefore $\mu \in \Omega^1_B ({\mathcal F})$.
Similarly (by $\eta (\xi ) = 0$)
\[ 2 (\xi \, \rfloor \, d \eta ) X = \xi (\eta (X)) - X(\eta (\xi
)) - \eta ([\xi , X]) = \] \[ = \xi (g(X,W)) - g({\mathcal L}_\xi
X , W) = ({\mathcal L}_\xi g)(X,W) + g(X, {\mathcal L}_\xi W ) = 0
\]
so that $\eta \in \Omega^1_B ({\mathcal F})$. Finally the
identities
\[ \xi \, \rfloor \, (d f \; \wedge \mu \wedge \eta ) =
\frac{1}{3} \; \xi (f) \; \mu \wedge \eta , \]
\[ \xi \, \rfloor \, (d \mu \wedge \eta ) = \frac{2}{3} \; (\xi \,
\rfloor \, d \mu ) \wedge \eta = 0, \;\;\; \xi \, \rfloor \, \mu
\wedge d \eta = \frac{2}{3} \; \mu \wedge (\xi \, \rfloor \, d
\eta ) = 0, \] together with $\xi \, \rfloor \, d \Omega = 0$
yield $\xi (f) = 0$.
\par
Next note that
\begin{equation} d \mu = 0, \;\;\; d \eta = 0.
\label{e:fol20}
\end{equation}
Indeed $\xi \, \rfloor \, d \mu = 0$ yields $d \mu = h \; \mu
\wedge \eta$ for some $h \in C^\infty (M)$. On the other hand
\[ h = 2 (d \mu )(N, W) = \] \[ = N(\mu (W)) - W (\mu (N)) - \mu ([N,W]) =
\] \[ =
g(\xi , [W, N]) = 0. \] The proof that $d \eta = 0$ is similar.
Let now $\omega \in \Omega^1_B ({\mathcal F})$. Then $\xi \,
\rfloor \, \omega = 0$ yields $\omega = a \mu + b \eta$ for some
$a,b \in C^\infty (M)$. By (\ref{e:fol20})
\[ d \omega = d a \; \wedge \mu + d b \; \wedge \eta = \]
\[ = (\xi (a) \lambda + W(a) \eta ) \wedge \mu + (\xi (b) \lambda
+ N(b) \mu ) \wedge \eta = \]
\[ = \xi (a) \, \lambda \wedge \mu + \xi (b) \, \lambda \wedge
\eta + (N(b) - W(a)) \, \mu \wedge \eta \] and $0 = 2 \; \xi \,
\rfloor \, d \omega = \xi (a) \, \mu + \xi (b) \, \eta$ shows that
both $a$ and $b$ are basic functions. Finally
\[ H^2_B ({\mathcal F}) = \frac{{\rm Ker}(d : \Omega^2_B
({\mathcal F}) \to \cdot )}{d \Omega^1_B ({\mathcal F})} = \] \[ =
\frac{\{ f \; \mu \wedge \eta : f \in \Omega^0_B ({\mathcal
F})\}}{\{ (N(b) - W(a)) \, \mu \wedge \eta : a,b \in \Omega^0_B
({\mathcal F}) \}} = 0. \] Proposition \ref{p:fol3} is proved.
\begin{corollary} Let $(M , g)$ be a complete $3$-dimensional Lorentz manifold
with ${\rm Isom}(M , g) = {\rm O}(4,1)$. If $M \sim S^3_1 (r)$
i.e. $M$ is a real homology pseudosphere $S^3_1 (r)$ then $M$
admits no complemented lightlike Killing vector field.
\label{c:fol1}
\end{corollary}
Here $S^n_\nu (r)$ is the pseudosphere i.e. $S^n_\nu (r) = \{ x
\in {\mathbb R}^{n+1} : - \sum_{j=1}^\nu x_j^2 + \sum_{j=\nu +
1}^{n+1} x_j^2 = r^2 \}$ ($r > 0$). As well known \[ H^j (S^n_\nu
(r), {\mathbb R}) = \begin{cases} {\mathbb R} & {\rm if} \; j \in
\{ 0, n-\nu \} , \cr 0 & {\rm otherwise} . \rm \end{cases} \] The
proof of Corollary \ref{c:fol1} is by contradiction. Let $\xi$ be
a complemented lightlike Killing vector field on $M$ and let $H$
be the $1$-parameter group of transformations generated by $\xi$.
Its closure $G = \overline{H} \subset {\rm O}(4,1)$ is compact (by
Lemma 3.1 in \cite{kn:Kam}) hence $G$ is a torus. The proof of
Theorem \ref{t:fol1} relies only on the existence of $N \in
S(T{\mathcal F}^\bot )^\bot$ such that $g(\xi , N) = 1$, $g(N,N) =
0$ and ${\mathcal L}_\xi N = 0$ so that we obtain the long exact
cohomology sequence (\ref{e:fol17}). In particular
\begin{equation}  H^1 (M , {\mathbb R}) \to H^0_B
({\mathcal F}) \to H^2_B ({\mathcal F}) \to H^2 (M , {\mathbb R})
\to H^1_B ({\mathcal F}) \label{e:fol21}
\end{equation}
is exact. Since $M$ is assumed to have the real cohomology of
$S^3_1 (r)$ one has $H^0 (M , {\mathbb R}) = H^2 (M , {\mathbb R})
= {\mathbb R}$ and $H^1 (M , {\mathbb R}) = 0$. Thus $H^1_B
({\mathcal F}) = 0$ and (\ref{e:fol21}) yields the exact sequence
\[ 0 \to {\mathbb R} \to H^2_B ({\mathcal F}) \to {\mathbb R} \to 0, \]
in contradiction with Proposition \ref{p:fol3}. Corollary
\ref{c:fol1} is proved. Under the same assumptions as those of
Theorem \ref{t:fol1} (with $M$ not necessarily compact) one also
has
\begin{corollary} {\rm 1)} If $M \sim S^{\nu +1}_\nu (r)$ then for
any $\ell \geq 1$
\[ H^{2\ell}_B ({\mathcal F}) = \begin{cases} {\mathbb R} & {\rm
if} \; \Delta (1) \neq 0, \cr 0 & {\rm if} \; \Delta (1) = 0, \cr
\end{cases} \hspace{0.5cm} H^{2\ell -1}_B ({\mathcal F}) =
\begin{cases} {\mathbb R} & {\rm if} \; H^1_B ({\mathcal F}) \neq 0,
\cr 0 & {\rm if} \; H^1_B ({\mathcal F}) = 0, \cr \end{cases} \]
where $\Delta : {\mathbb R} \to H^2_B ({\mathcal F})$ is the map
$\Delta (c) = c [d \alpha ]$, $c \in {\mathbb R}$. {\rm 2)} If $M
\sim S^{\nu +2}_\nu (r)$ then for any $\ell \geq 0$
\[ H^{2\ell +1}_B ({\mathcal F}) = 0, \] \[ H^{2\ell}_B
({\mathcal F}) = \begin{cases} {\mathbb R} & {\rm if} \; \exists
\; f \in C^\infty (M) \; {\rm with} \; \xi (f) = 1, \cr 2-{\rm
dimensional} & {\rm otherwise}. \cr \end{cases} \] {\rm 3)} Assume
that $M \sim S^{\nu + p}_\nu (r)$ for some fixed $3 \leq p \leq
n-1$. Then
\[ H^j_B ({\mathcal F}) = \begin{cases} {\mathbb R} & {\rm if} \;
j = {\rm even}, \cr 0 & {\rm if} \; j = {\rm odd}, \cr \end{cases}
\;\;\; 0 \leq j \leq p-1, \]
\[ H^{p+2j}_B ({\mathcal F}) \approx H^p_B ({\mathcal F} ), \;\;\;
1 \leq j \leq \left[ \frac{n-p-1}{2} \right] , \]
\[ H^{p + 2j+1}_B ({\mathcal F}) \approx H^{p+1}_B ({\mathcal F}),
\;\;\; 1 \leq j \leq \left[\frac{n-p}{2}\right] - 1. \] Moreover
one has either that {\rm i)} $j_* H^p ({\mathcal F}) = 0$ and then
$H^{p-2}_B ({\mathcal F}) \approx H^p ({\mathcal F})$, or that
{\rm ii)} $j_* H^p_B ({\mathcal F}) = {\mathbb R}$ and then either
$H^{p-2}_B ({\mathcal F}) = 0$ and $H^p_B ({\mathcal F}) =
{\mathbb R}$ or $H^{p-2}_B ({\mathcal F}) \neq 0$ and
$\dim_{\mathbb R} H^p_B ({\mathcal F})  = 2$. Finally either {\rm
iii)} $(i_\xi )_* (1) = 0$ and then $H^{p+1}_B ({\mathcal F})
\approx H^{p-1}_B ({\mathcal F})$, or {\rm iv)} $(i_\xi )_* (1)
\neq 0$ and then $\dim_{\mathbb R} H^{p+1}_B ({\mathcal F}) =
\dim_{\mathbb R} H^{p-1}_B ({\mathcal F}) - 1$. \label{c:folhom}
\end{corollary}
Here $[a] \in {\mathbb Z}$ is the integer part of $a \in {\mathbb
R}$. A similar result holds when $M$ is a real homology
pseudohyperbolic space. {\em Proof of Corollary} \ref{c:folhom}.
Statement 1) in Corollary \ref{c:folhom} is a straightforward
consequence of (\ref{e:fol17}) and $H^j (M) = {\mathbb R}$ if $n =
\nu + j$ ($1 \leq j \leq n-1$) and $H^j (M) = 0$ otherwise.
Statement 2) requires again (\ref{e:fol17}) and the following
\begin{lemma} If there is $f \in C^\infty (M)$ such that $\xi (f)
= 1$ then $H^2_B ({\mathcal F}) = {\mathbb R}$. Otherwise
$\dim_{\mathbb R} H^2_B ({\mathcal F}) = 2$. \label{l:folhom}
\end{lemma}
{\em Proof}. Let $f \in C^\infty (M)$ such that $\xi (f) = 1$ and
let us set $\omega := \alpha - d f$. Then $\xi \, \rfloor \,
\omega = 0$. Also $d \omega = d \alpha$ so that $\xi \, \rfloor \,
d \omega = 0$. Hence $\omega \in \Omega^1_B ({\mathcal F})$ i.e.
$\omega$ is a basic form and then $\Delta (1) = [d \alpha ] = [d
\omega ] = 0 \in H^2_B ({\mathcal F})$. Note that $H^1_B
({\mathcal F}) \hookrightarrow H^1 (M) = 0$. Then $\Delta
({\mathbb R}) = 0$ and the exactness of
\[ {\mathbb R} \stackrel{\Delta}{\to} H^2_B ({\mathcal F})
\stackrel{j_*}{\to} {\mathbb R} \to 0 \] imply that $j_* : H^2_B
({\mathcal F}) \approx {\mathbb R}$. Otherwise $\Delta ({\mathbb
R}) \neq 0$ and $H^2_B ({\mathcal F})/\Delta ({\mathbb R}) \approx
{\mathbb R}$ so that Lemma \ref{l:folhom} is proved. Finally,
statement 3) in Corollary \ref{c:folhom} is implied by
(\ref{e:fol17}) and the fact that all cohomology groups of $M$
vanish except for $H^p (M) \approx {\mathbb R}$.

\section{Rummler's formula for lightlike foliations
with trivial radical distribution} Let $({\mathcal F}, \;
S(T{\mathcal F}), \; S(T{\mathcal F}^\bot ))$ be a $r$-lightlike
foliation of the semi-Rieman\-nian manifold $M$, together with a
choice of screen distributions. We say that $\mathcal F$ is {\em
tangentially } (respectively {\em transversally}) {\em screen
oriented} if $S(T{\mathcal F})$ is oriented (respectively if
$S(T{\mathcal F}^\bot )$ is oriented). As both screen
distributions are nondegenerate we may consider the local
orthonormal frames $\{ X_a : 1 \leq a \leq m-r \}$ in
$S(T{\mathcal F})$ and $\{ W_\alpha : 1 \leq \alpha \leq q - r \}$
in $S(T{\mathcal F}^\bot )$, defined on the open set $U \subseteq
M$. Then $g(X_a , X_b ) = \epsilon_a \delta_{ab}$ and $g(W_\alpha
, W_\beta ) = \epsilon_\alpha \delta_{\alpha\beta}$ (where
$\epsilon_j^2 = 1$). Moreover we set $\omega^a (X) = g(X , X_a )$
and $\eta^\alpha (X) = g(X , W_\alpha )$ for any $X \in T(M)$. If
$\mathcal F$ is tangentially screen oriented then the local
$(m-r)$-forms $\omega^1 \wedge \cdots \wedge \omega^{m-r}$ glue up
to give a (globally defined) $(m-r)$-form $\chi_{S(T{\mathcal
F})}$ on $M$.
\par
We shall need a lightlike counterpart of {\em Rummler's formula},
cf. H. Rummler, \cite{kn:Rum} (or identity (6.17) in
\cite{kn:Ton}, p. 66). To this end one ought to build a lightlike
analog of the characteristic form $\chi_{\mathcal F}$ of a
tangentially oriented foliation $\mathcal F$ of a Riemannian
manifold, cf. \cite{kn:Ton}, p. 65-66. We do this under the
additional assumption that
\[ {\rm Rad} \; T{\mathcal F} \approx M \times {\mathbb R}^r \]
(the trivial bundle). If this is the case,
let $\{ \xi_i : 1 \leq i \leq r \}$ be a globally defined frame
(fixed through the remainder of this section) of ${\rm Rad} \;
T{\mathcal F}$ and let $\{ N_i : 1 \leq i \leq r \}$ be the
lightlike vector fields furnished by Lemma \ref{l:fol1}. The
construction of $N_i$ depends on a choice of complement $E$ to
${\rm Rad} \; T{\mathcal F}$ in $S(T{\mathcal F}^\bot )^\bot$ as
in (\ref{e:fol7}) and on the choice of a {\em local} frame $\{ V_i
: 1 \leq i \leq r \}$ of $E$ on $U \subseteq M$, so that {\em a
priori} the $N_i$'s are but locally defined. Nevertheless if $E$
is fixed as well and $\{ V^\prime_i : 1 \leq i \leq r \}$ is
another local frame of $E$ on $U^\prime \subseteq M$ then
$N^\prime_i = N_i$ on $U \cap U^\prime$. Therefore the proof of
Lemma \ref{l:fol1} yields a {\em globally defined} system of
vector fields $N_i \in S(T{\mathcal F}^\bot )^\bot$ such that
$g(\xi_i , N_j ) = \delta_{ij}$ and $g(N_i , N_j ) = 0$. We set
$\lambda^i (X) = g(X , N_i )$ and $\mu^i (X) = g(X, \xi_i )$ for
any $X \in T(M)$. Moreover let $\chi_{\mathcal F}$ be the $m$-form
on $M$ given by
\[ \chi_{\mathcal F} = \lambda^1 \wedge \cdots \wedge
\lambda^r \wedge \chi_{S(T{\mathcal F})}. \] We emphasize on a
number of elementary properties of the local frame $\{ X_a , \;
\xi_i , W_\alpha , N_i \}$. First $g(N_i , N_j ) = 0$ may be
written
\begin{equation} \lambda^i (N_j ) = 0, \;\;\; 1 \leq i,j \leq r.
\label{e:tr1}
\end{equation}
Next
\[ W_\beta \in S(T {\mathcal F}^\bot ) \subset T({\mathcal F})^\bot \; \bot \;
T({\mathcal F}) \supset S(T{\mathcal F}) \ni X_a \] hence $g(X_a ,
W_\beta ) = 0$ i.e.
\begin{equation} \omega^a (W_\beta ) = 0, \;\;\; 1 \leq a \leq
m-r, \;\;\; 1 \leq \beta \leq q - r. \label{e:tr2}
\end{equation}
Moreover
\[ N_i \in ltr(T{\mathcal F}) \subset S(T{\mathcal F}^\bot )^\bot
\; \bot \; S(T{\mathcal F}^\bot ) \ni W_\beta \] hence $g(N_i ,
W_\beta ) = 0$ i.e.
\begin{equation}
\lambda^i (W_\beta ) = 0, \;\;\; 1 \leq i \leq r, \;\;\; 1 \leq
\beta \leq q - r. \label{e:tr3}
\end{equation}
Also
\[ N_i \in ltr(T {\mathcal F}) \subset S(T{\mathcal F})^\bot \;
\bot \; S(T{\mathcal F}) \ni X_a \] so that $g(X_a , N_j ) = 0$
i.e.
\begin{equation} \omega^a (N_j ) = 0, \;\;\; 1 \leq a \leq m-r,
\;\;\; 1 \leq j \leq r. \label{e:tr4}
\end{equation}
Then (\ref{e:tr1})-(\ref{e:tr4}) imply that
\begin{equation}
tr(T{\mathcal F}) \, \rfloor \, \lambda^i = 0, \;\;\;
tr(T{\mathcal F}) \, \rfloor \, \omega^a = 0.
\label{e:tr5}
\end{equation}
In particular $tr(T{\mathcal F}) \, \rfloor \, \chi_{\mathcal F} =
0$ and
\[ \chi_{\mathcal F}(\xi_1 , \cdots ,
\xi_r , X_1 , \cdots , X_{m-r}) = 1/m! \, . \]
Let ${\rm tan} :
T(M) \to T({\mathcal F})$ be the natural projection associated
with the decomposition (\ref{e:fol11}). If $Z \in \Gamma^\infty
(tr(T{\mathcal F}))$ then
\[ ({\mathcal L}_Z \chi_{\mathcal F} )(Y_1 , \cdots , Y_m ) = Z
(\chi_{\mathcal F}(Y_1 , \cdots , Y_m )) - \]
\[ - \sum_{j=1}^m \chi_{\mathcal F} (Y_1 , \cdots , {\rm tan} [Z,
Y_j ] , \cdots , Y_m ) \] for any $Y_j \in T({\mathcal F})$. We
wish to evaluate this identity at $Y_i = \xi_i$ and $Y_{\alpha +
r} = X_\alpha$. The first term on the right hand side vanishes. By
(\ref{e:tr5}) one has ${\rm tan}(X) = \lambda^i (X) \xi_i +
\epsilon^a \omega^a (X) X_a$ for any $X \in T(M)$ (where
$\epsilon^a = \epsilon_a$). Hence
\begin{proposition} Let $({\mathcal F}, \; S(T{\mathcal F}), \;
S(T{\mathcal F}^\bot ))$ be a $r$-lightlike foliation of
codimension $q$ of a $n$-dimensional semi-Riemannian manifold $(M
, g)$ such that $1 \leq r \leq m-1$ and $1 \leq r \leq q-1$ where
$m = n-q$. Let us assume that ${\rm Rad} \; T{\mathcal F}$ is
trivial and let $\xi = (\xi_1 , \cdots ,  \xi_r )$ be a global
frame of ${\rm Rad} \; T{\mathcal F}$. If we define $\kappa =
\kappa (\xi , E)$ by setting
\[ \kappa (X) = 0, \;\;\; X \in
T({\mathcal F}), \]
\[ \kappa (Z) =
\lambda^i ([Z, \xi_i ]) + \epsilon^a \omega^a ([Z, X_a ]), \;\;\;
Z \in tr(T{\mathcal F}), \]  then $\kappa \in \Omega^1 (M)$ i.e.
$\kappa$ is globally defined and its definition doesn't depend
upon the choice of local orthonormal frames $\{ X_a \} \subset
S(T{\mathcal F})$ and $\{ W_\alpha \} \subset S(T{\mathcal F}^\bot
)$. If $\mathcal F$ is tangentially screen oriented then
\begin{equation}
{\mathcal L}_Z \chi_{\mathcal F} + \kappa (Z) \chi_{\mathcal F} =
0 \label{e:fol22}
\end{equation}
on $T({\mathcal F}) \otimes \cdots \otimes T({\mathcal F})$ for
any $Z \in tr(T{\mathcal F})$. \label{p:kappa}
\end{proposition}
The identity (\ref{e:fol22}) is the lightlike analog to the
Rummler formula we were seeking for while $\kappa$ is formally
similar to the mean curvature form of a Riemannian foliation (cf.
e.g. \cite{kn:Ton}, p. 67). Let $h$ be the second fundamental form
of the foliation $\mathcal F$ i.e.
\[ h (X,Y) = \Pi \nabla^g_X Y , \;\;\; X,Y \in T({\mathcal
F}). \] See also (6.1) in \cite{kn:Ton}, p. 62. Let
$h_{S(T{\mathcal F})}$ be the restriction of $h$ to $S(T{\mathcal
F}) \otimes S(T{\mathcal F})$ and $\tau_{S(T{\mathcal F})}$ the
trace of $h_{S(T{\mathcal F})}$ with respect to $g$ i.e.
\[ \tau_{S(T{\mathcal F})} = {\rm trace}_g \; h_{S(T{\mathcal F})}
\] (locally $\tau_{S(T{\mathcal F})} = \sum_{a=1}^{m-r} \epsilon_a
h(X_a , X_a )$). As it turns out, in the case of lightlike
foliations of semi-Riemannian manifolds neither
$\tau_{S(T{\mathcal F})}$ is the mean curvature vector of the
distribution $S(T{\mathcal F})$ in $(M , g)$ nor $\kappa$ equals
$\tau_{S(T{\mathcal F})}$ (but rather $\tau_{S(T{\mathcal F})}$
plus extra terms whose geometric meaning is rather obscure).
Indeed (as $\nabla^g$ is torsion-free and $\nabla^g g = 0$)
\[ \kappa (Z) = \lambda^i ([Z, \xi_i ]) + \epsilon^a g(\nabla^g_Z X_a , X_a ) - \epsilon^a
g(\nabla_{X_a} Z , X_a ) = \]
\[ = \lambda^i ([Z, \xi_i ]) +
\epsilon^a g(Z , \nabla^g_{X_a} X_a )  \] that is
\begin{equation}
\label{e:fol23} \kappa (Z) = \lambda^i ([Z, \xi_i ]) + g (Z ,
H_{S(T{\mathcal F})})
\end{equation}
where $H_{S(T{\mathcal F})}$ is the mean curvature vector of
$S(T{\mathcal F})$ in $(M , g)$ i.e.
\[ H_{S(T{\mathcal F})} = {\rm trace}_g B_{S(T{\mathcal F})} \, ,
\]
\[ B_{S(T{\mathcal F})} (X,Y) = (\nabla^g_X Y )_{S(T{\mathcal
F})^\bot} \, , \;\;\; X,Y \in S(T{\mathcal F}), \] and
$V_{S(T{\mathcal F})^\bot}$ is the $S(T{\mathcal
F})^\bot$-component of $V \in T(M)$ with respect to the direct sum
decomposition  (\ref{e:fol3}). Indeed, as
\[ S(T{\mathcal F})^\bot = tr(T{\mathcal F}) \oplus {\rm Rad} \;
T{\mathcal F} \] $S(T{\mathcal F})^\bot$ is locally the span of
$\{ N_i , W_\alpha , \xi_i \}$ hence \[ g(H_{S(T{\mathcal F})} ,
Z) = \sum_{a=1}^{m-r} \epsilon^a g(\nabla^g_{X_a} X_a \, , \, Z)
\]
and (\ref{e:fol23}) is proved. The identity (\ref{e:fol23}) also
shows that $\kappa$ is indeed globally defined. On the other hand,
if $s := \Pi (Z)$ then
\[ g_{\rm tra} (s , \tau_{S(T{\mathcal F})}) = \sum_a \epsilon^a
g(Z , {\rm tra} (\nabla^g_{X_a} X_a )) \] so that
\[ g(Z, H_{S(T{\mathcal F})}) = g_{\rm tra}(s, \tau_{S(T{\mathcal
F})} ) + \epsilon^a g({\rm ltr}(Z) , \nabla^g_{X_a} X_a ) \] where
${\rm ltr} : tr(T{\mathcal F}) \to ltr(T{\mathcal F})$ is the
projection associated to the decomposition (\ref{e:fol10}).
\par
We shall need the multiplicative filtration $\{ F^r \Omega^k : r
\geq 0, \; k \geq r - 1 \}$ of the de Rham complex $\Omega^\bullet
(M)$ as devised by F. Kamber \& P. Tondeur, \cite{kn:KaTo} (cf.
also \cite{kn:Ton}, p. 120). That is
\[ F^r \Omega^k = \{ \omega \in \Omega^k (M) : i_{X_1} \cdots
i_{X_{k-r+1}} \omega = 0, \;\;\; X_j \in T({\mathcal F}) \} . \]
An useful reformulation of (\ref{e:fol22}) is ${\mathcal L}_Z
\chi_{\mathcal F} + \kappa (Z) \chi_{\mathcal F} \in F^1 \Omega^m$
or equivalently (as $T({\mathcal F}) \, \rfloor \, \kappa = 0$)
\begin{equation}  d \chi_{\mathcal F} + (m+1) \kappa \wedge
\chi_{\mathcal F} \in F^2 \Omega^{m+1} \, .
\label{e:fol24}
\end{equation}
The way $\kappa$ depends upon $(\xi , E)$ is described by the
following
\begin{corollary} Let $\mathcal F$ be a $r$-lightlike foliation
of a semi-Riemannian manifold with ${\rm Rad} \; T {\mathcal F}$
trivial. Let $\xi^\prime = (\xi^\prime_1 , \cdots , \xi^\prime_r
)$ be another global frame of ${\rm Rad} \; T{\mathcal F}$ and let
$E^\prime \to M$ be another complement to ${\rm Rad} \; T
{\mathcal F}$ in $S(T{\mathcal F}^\bot )^\bot$. Then there exists
a $C^\infty$ function $f : M \to {\rm GL}(r, {\mathbb R})$ such
that
\begin{equation}
\kappa (\xi^\prime , E^\prime ) = \kappa (\xi , E) + {\rm trace}
(f^{-1} \, d f). \label{e:kappa}
\end{equation}
In particular if $\mathcal F$ is a $1$-lightlike foliation and
there is $(\xi , E)$ such that $d \; \kappa (\xi , E) = 0$ then
the de Rham cohomology class $[\kappa (\xi , E)] \in H^1 (M ,
{\mathbb R})$ doesn't depend upon the choice of $(\xi , E)$. For
instance let $M$ be complete and let $\xi$ be a lightlike Killing
vector field on $M$. Let $\mathcal F$ be the $1$-lightlike
foliation of $M$ such that $T({\mathcal F}) = {\mathbb R} \xi$.
Assume that there exist $G$-invariant globally defined nowhere
zero vector fields $V \in S(T{\mathcal F}^\bot )^\bot$ and $W \in
S(T{\mathcal F}^\bot )$ such that {\rm (\ref{e:fol16})} holds.
Then $\kappa (\xi , {\mathbb R}V )$ is closed.
\end{corollary}
{\em Proof}. Let $\{ V^\prime_i : 1 \leq i \leq r \}$ be a local
frame of $E^\prime$, defined on an open set $U^\prime \subseteq M$
such that $U \cap U^\prime \neq \emptyset$. Then $V^\prime_i =
a_i^j \xi_j + b_i^j V_j$ for some $C^\infty$ functions $a_i^j , \,
b_i^j : U \cap U^\prime \to {\mathbb R}$. Let us set
$g_{ij}^\prime = g(\xi_i^\prime , V_j^\prime )$ and (according to
the proof of Lemma \ref{l:fol1}) $[{g^\prime}^{ij}] =
[g_{ij}^\prime ]^{-1}$. Let $f = [f^i_j ] : M \to {\rm GL}(r ,
{\mathbb R})$ such that $\xi^\prime = f \xi$. Then
\begin{equation}
g^{ij} = b^i_\ell \, f^j_k \, {g^\prime}^{k\ell} . \label{e:gij}
\end{equation}
Let $N^\prime_j$ be given by
\[ N^\prime_j = - \frac{1}{2} {g^\prime}^{ki} {g^\prime}^{\ell j}
g(V^\prime_k , V^\prime_\ell ) \xi_i^\prime + {g^\prime}^{kj}
V^\prime_k \, , \;\;\; 1 \leq j \leq r.  \] A calculation based on
(\ref{e:gij}) shows that
\begin{equation}
N^\prime_j = \sum_{i=1}^r (f^{-1})^j_i N_i + \frac{1}{2} \{ a^i_k
{g^\prime}^{kj} - a_k^p {g^\prime}^{k\ell} (f^{-1})^j_p f^i_\ell
\} \xi_i \, . \label{e:Nprime}
\end{equation}
Consequently
\[ {\lambda^\prime}^j ([Z , \xi^\prime_j ]) = \lambda^i ([Z, \xi_i
]) + (f^{-1})^j_i \, Z(f^i_j ) \] yielding (\ref{e:kappa}). When
$r = 1$ the identity (\ref{e:kappa}) becomes \[ \kappa (\xi^\prime
, E^\prime ) = \kappa (\xi , E) + d \log |f|. \]
Q.e.d.

\section{The transversal divergence theorem}
Assume from now on that $\mathcal F$ is transversally screen
oriented. Then the local $(q-r)$-forms $\eta^1 \wedge \cdots
\wedge \eta^{q-r}$ glue up to a (globally defined) $(q-r)$-form
$\nu_{S(T{\mathcal F}^\bot )}$ on $M$. Let $\nu_{\mathcal F}$ be
the $q$-form on $M$ given by $\nu_{\mathcal F} = \mu^1 \wedge
\cdots \wedge \mu^r \wedge \nu_{S(T{\mathcal F}^\bot )}$. Then
$\omega := \nu_{\mathcal F} \wedge \chi_{\mathcal F}$ is a volume
form on $M$. We denote by $V({\mathcal F})$ the set of all
infinitesimal automorphisms of $\mathcal F$ i.e. $Y \in
V({\mathcal F})$ is a vector field on $M$ such that $[X,Y] \in
T({\mathcal F})$ for any $X \in T({\mathcal F})$. The transversal
divergence operator is the map ${\rm div}_B : V({\mathcal F}) \to
C^\infty (M)$ given by
\[ {\mathcal L}_Y \nu_{\mathcal F} = {\rm div}_B (Y) \;
\nu_{\mathcal F} , \;\;\; Y \in V({\mathcal F}). \] Cf. e.g.
\cite{kn:Ton}, p. 126. One checks easily that
\begin{lemma} Assume that {\rm 1)}
the complement $E$ in {\rm (\ref{e:fol7})} and the local frame $\{
V_i : 1 \leq i \leq r \}$ may be chosen such that $N_i \in
V({\mathcal F})$, $1 \leq i \leq r$, and {\rm 2)} there is a
transversal screen distribution $S(T{\mathcal F}^\bot )$ admitting
a local orthonormal frame $\{ W_\alpha : 1 \leq \alpha \leq q- r
\}$ such that each $W_\alpha$ is a local infinitesimal
automorphism of $\mathcal F$. Then $\nu_{\mathcal F}$ is holonomy
invariant. In particular ${\rm div}_B X = 0$ for any $X \in
T({\mathcal F})$ and ${\rm div}_B Y \in \Omega^0_B ({\mathcal F})$
for any $Y \in V({\mathcal F})$. \label{l:fol3}
\end{lemma}
We shall prove the following lightlike analog of a result by F.
Kamber \& P. Tondeur \& G. Toth, \cite{kn:KTT}
\begin{theorem} Let ${\mathcal F}$ be a tangentially and transversally screen
oriented lightlike foliation of the semi-Riemannian manifold $(M,
g)$. Assume that the radical distribution of $\mathcal F$ is
trivial {\rm (}i.e. ${\rm Rad} \; T{\mathcal F} \approx M \times
{\mathbb R})$ and there is a transverse bundle $tr(T{\mathcal F})
= ltr(T{\mathcal F}) \oplus S(T{\mathcal F}^\bot )$ with
$ltr(T{\mathcal F})$ and $S(T{\mathcal F}^\bot )$ as in Lemma {\rm
\ref{l:fol3}}. If $\partial M = \emptyset$ then
\begin{equation} \int_M {\rm div}_B (Y) \; \nu_{\mathcal F} \wedge
\chi_{\mathcal F} = (-1)^q (m+1) \int_M (i_{Y} \kappa ) \;\;
\nu_{\mathcal F} \wedge \chi_{\mathcal F}
\end{equation}
for any compactly supported $Y \in V({\mathcal F})$.
\label{t:fol2}
\end{theorem}
When $M$ is a manifold-with-boundary ($\partial M \neq \emptyset$)
the problem of producing a foliated analog to the semi-Riemannian
divergence theorem (cf. B. \"Unal, \cite{kn:Una}) is left open.
\par
{\em Proof of Theorem} \ref{t:fol2}. Let $Y \in V({\mathcal F})$
and $Z := {\rm tra}(Y)$. A calculation shows that $(i_Z \; d \;
\nu_{\mathcal F}) \wedge \chi_{\mathcal F} = 0$. Then we have (by
Lemma \ref{l:fol3})
\[ ({\rm div}_B Y ) \omega =
({\rm div}_B Z) \nu_{\mathcal F} \wedge \chi_{\mathcal F} =
({\mathcal L}_Z \nu_{\mathcal F}) \wedge \chi_{\mathcal F} =
\]
\[ = \{ (d i_Z + i_Z d) \nu_{\mathcal F} \} \wedge \chi_{\mathcal
F} = (d \; i_Z \; \nu_{\mathcal F}) \wedge \chi_{\mathcal F} = \]
\[ = d \left( (i_Z \nu_{\mathcal F} ) \wedge \chi_{\mathcal F}
\right) + (-1)^q (i_Z \nu_{\mathcal F}) \wedge d \chi_{\mathcal F}
. \] If we set $\varphi := d \chi_{\mathcal F} + (m+1) \kappa
\wedge \chi_{\mathcal F}$ then (as $T({\mathcal F}) \, \rfloor \,
\nu_{\mathcal F} = 0$)
\begin{equation}
\label{e:fol27} ({\rm div}_B Y ) \omega = d \left( (i_Y
\nu_{\mathcal F} ) \wedge \chi_{\mathcal F} \right) + (-1)^q (i_Y
\nu_{\mathcal F}) \wedge (\varphi - (m+1) \kappa \wedge
\chi_{\mathcal F}) .
\end{equation}
Note that $\varphi \in F^2
\Omega^{m+1}$ (by (\ref{e:fol24})). Let $(U, x^1 , \cdots , x^m ,
y^1 , \cdots , y^q )$ be a foliated chart. Since $T({\mathcal F})
\, \rfloor \, i_Y \nu_{\mathcal F} = 0$ it follows that $i_Y
\nu_{\mathcal F}$ is a sum of monomials of the form $d y^1 \wedge
\cdots \wedge d y^{q-1}$ with $C^\infty$ coefficients. Also
$\varphi$ is a sum of monomials each of which contains at least a
monomial of the form $d y^i \wedge d y^j$. Hence
\[ (i_Y \nu_{\mathcal F} ) \wedge \varphi \in F^{q+2} \Omega^n =
0 \] and (\ref{e:fol27}) becomes
\[ ({\rm div}_B Y) \omega = d
\left( (i_Y \nu_{\mathcal F} ) \wedge \chi_{\mathcal F} \right) +
(-1)^{q+1} (m+1) (i_Y \nu_{\mathcal F}) \wedge \kappa \wedge
\chi_{\mathcal F} . \] Next $T({\mathcal F}) \, \rfloor \, \kappa
= 0$ yields $\kappa \wedge \nu_{\mathcal F} = 0$ and then $(i_Y
\kappa ) \wedge \nu_{\mathcal F} - \kappa \wedge i_Y \nu_{\mathcal
F} = 0$ so that
\begin{equation} ({\rm div}_B Y) \omega =
d \left( (i_Y \nu_{\mathcal F}) \wedge \chi_{\mathcal F} \right) +
(-1)^q (m+1)(i_Y \kappa ) \omega . \label{e:fol28}
\end{equation}
Finally one integrates (\ref{e:fol28}) over $M$ and uses the
Stokes theorem.

\section{Lightlike functions}
Let $(M , g)$ be a $n$-dimensional semi-Riemannian manifold. A
$C^\infty$ function $f : M \to {\mathbb R}$ is said to be {\em
lightlike} if $\nabla f$ is null i.e. ${\rm Crit}(f) = \emptyset$
and $g(\nabla f , \nabla f) = 0$. For instance if $M = {\mathbb
R}^n_s$ ($1 \leq s \leq n-1$) then a smooth function $f : {\mathbb
R}^n_s \to {\mathbb R}$ is lightlike if $(f_{x_1}(x) , \cdots ,
f_{x_n}(x)) \neq 0$ at any $x \in {\mathbb R}^n$ and \[
\sum_{j=1}^s \left( f_{x_j} \right)^2 - \sum_{j=s+1}^n \left(
f_{x_j} \right)^2 = 0, \] where $f_{x_j} =
\partial f/\partial x_j$. Let $f : M \to {\mathbb R}$ be a
lightlike function and let $\mathcal F$ be the foliation by level
sets of $f$ so that \[ T({\mathcal F}) = \{ X \in T(M) : X(f) = 0
\} . \]
Then $\dim_{\mathbb R} T({\mathcal F})^\bot_x = 1$.
Consequently $T({\mathcal F})^\bot = {\mathbb R} \nabla f$. Taking
into account the classification in Table 1
\begin{proposition} For any lightlike function on a
semi-Riemannian manifold the corresponding foliation by level sets
is co-isotropic.
\end{proposition}
\noindent  Indeed ${\rm Rad} \; T {\mathcal F} = T({\mathcal F})
\cap T({\mathcal F})^\bot = {\mathbb R} \nabla f$. We wish to
apply the results in Section 1 to foliations by level sets of
lightlike functions. Therefore we choose a screen distribution
$S(T{\mathcal F})$ such that
\begin{equation} T({\mathcal F}) =
S(T{\mathcal F}) \oplus {\mathbb R} \nabla f \label{e:fol39}
\end{equation}
(and $S(T{\mathcal F}^\bot ) = 0$). As $S(T{\mathcal F})$ is
nondegenerate $T(M) = S(T{\mathcal F}) \oplus S(T{\mathcal
F})^\bot$. In particular $\dim_{\mathbb R} S(T{\mathcal F})^\bot =
2$. Let $E \to M$ be a real line bundle such that
\begin{equation}
S(T{\mathcal F})^\bot = ({\mathbb R} \nabla f) \oplus E.
\label{e:fol40}
\end{equation}
Given $V \in \Gamma^\infty (U, E)$ such that $V_x \neq 0$ for any
$x \in M$ we set
\begin{equation}
N = \frac{1}{V(f)} \, \{ V - \frac{1}{2} \; \frac{g(V,V)}{V(f)} \;
\nabla f \} \label{e:fol41}
\end{equation}
hence
\begin{equation} g(N,N) =0, \;\;\; N(f) = 1. \label{e:fol31}
\end{equation}
Clearly if $N^\prime$ is similarly built in terms of a nowhere
zero $V^\prime \in \Gamma^\infty (U^\prime , E)$ then $V^\prime =
\lambda V$ for some $C^\infty$ function $\lambda : U \cap U^\prime
\to {\mathbb R}$ hence $N = N^\prime$ on $U \cap U^\prime$. This
furnishes a globally defined $C^\infty$ section $N$ in
$S(T{\mathcal F})^\bot$ possessing the properties (\ref{e:fol31}).
Then $tr(T{\mathcal F}) = {\mathbb R} N$ is a choice of
transversal bundle and in particular
\[ T(M) = T({\mathcal F}) \oplus {\mathbb R} N. \]
\begin{proposition} Let $f : M \to {\mathbb R}$ be a lightlike function on the semi-Riemannian
manifold $(M , g)$. Let $\mathcal F$ be the foliation of $M$ by
level sets of $f$. Then the isomorphism $\sigma : \nu ({\mathcal
F}) \approx tr(T{\mathcal F})$ is given by $\sigma (s) = Y(f) N$
for any $s = \Pi (Y)$, $Y \in T(M)$. Consequently the second
fundamental form $h$ of $\mathcal F$ in $(M , g)$ is given by
\begin{equation}
\sigma \; h(X,Y) = - {\rm Hess}_f (X,Y) N, \;\;\; X,Y \in
T({\mathcal F}).
\end{equation}
Finally $\kappa \in \Omega^1 (M)$ is given by $T({\mathcal F}) \,
\rfloor \, \kappa = 0$ and
\begin{equation}
2 \kappa (N) = ({\mathcal L}_\xi g)(N , N) - \label{e:fol44}
\end{equation}
\[ - \frac{1}{V(f)} \{
{\rm trace}_g ({\mathcal L}_V g)_{S(T{\mathcal F})} -
\frac{g(V,V)}{V(f)} [ \square f - 2 {\rm Hess}_f (\xi , N)] \} \]
on $U \subseteq M$, where $\square$ is the Laplace-Beltrami
operator of $(M , g)$, while $\xi$ and ${\rm Hess}_f$ are the
gradient and Hessian of $f$. \label{p:fol8}
\end{proposition}
For instance
\begin{corollary}
Let $f : {\mathbb R}^n_s \to {\mathbb R}$ be the linear function
\begin{equation} f(x_1 , \cdots , x_n ) =
\sqrt{\frac{n-s}{s}} \sum_{i=1}^s x_i + \sum_{j=s+1}^n x_j
\label{e:fol45}
\end{equation}
and let $\mathcal F$ be the corresponding foliation of ${\mathbb
R}^n_s$ by affine hyperplanes. Let $S(T{\mathcal F})$ be the span
of $\{ X_a : 1 \leq a \leq n-2 \}$ where \begin{equation} X_a =
\begin{cases}
\partial_a - \sqrt{\frac{n-s}{s}} \,
\partial_n \, , & {\rm if} \;\; 1 \leq a \leq s-1, \cr
\partial_{a+1} - \partial_n \, , & {\rm if} \;\; s \leq a \leq
n-2. \cr \end{cases}
\label{e:fol46}
\end{equation}
Then $S(T{\mathcal F})$ is a screen distribution i.e. {\rm
(\ref{e:fol39})} holds. Also if
\begin{equation}
V = - \sqrt{\frac{n-s}{s}} \, \sum_{i=1}^{s-1} \partial_i +
\sum_{j=s}^n
\partial_j \label{e:fol47}
\end{equation}
then $E = {\mathbb R}V$ is a complement to ${\rm Rad} \;
T{\mathcal F}$ in $S(T{\mathcal F})^\bot$. Finally $V$ and $\nabla
f$ are Killing vector fields on ${\mathbb R}^n_s$ and consequently
$h = 0$ and $\kappa = 0$.
\end{corollary}
Here $\partial_i$ is short for $\partial /\partial x_i$. {\em
Proof of Proposition} \ref{p:fol8}. It suffices to compute $\kappa
(N)$. On one hand
\[ g([N, \xi ], N) = \frac{1}{2} \, ({\mathcal L}_\xi g)(N,N). \]
On the other hand
\[ \square f = {\rm trace}_g \; {\rm Hess}_f \, , \]
\[ {\rm Hess}_f (X,Y) = X(Y(f)) - (\nabla^g_X Y)(f), \;\;\; X,Y
\in T(M), \] and
\begin{lemma} Let $(M , g)$ be a $n$-dimensional semi-Riemannian manifold of index
$1 \leq s \leq n-1$. Let $\{ X_a : 1 \leq a \leq n-2 \}$ be a
local orthonormal {\rm (}i.e. $g(X_a , X_b ) = \epsilon_a
\delta_{ab}$, $\epsilon_a^2 = 1${\rm )} frame of $S(T{\mathcal
F})$ defined on the open set $U \subseteq M$. If $f : M \to
{\mathbb R}$ is a lightlike function and $\epsilon \in \{ \pm 1
\}$ then $\{ X_1 , \cdots , X_{n-1} , \xi \pm (\epsilon /2) \, N
\}$ is a local orthonormal frame of $T(M)$ on $U$. In particular a
screen distribution $S(T{\mathcal F})$ has index ${\rm
ind}(S(T{\mathcal F}), g) = s-1$.
\end{lemma}
Consequently
\[ - \epsilon^a (\nabla^g_{X_a} X_a )(f) = \square f - 2 \, {\rm
Hess}_f (\xi , N) \] hence (by (\ref{e:fol41}) and $V \in
S(T{\mathcal F})^\bot$)
\[ g(N , H_{S(T{\mathcal F})} ) = \frac{\epsilon^a}{V(f)} \{ g(V ,
\nabla^g_{X_a} X_a )- \frac{g(V,V)}{2 V(f)} \, (\nabla^g_{X_a} X_a
)(f) \} = \]
\[ = - \frac{1}{2 V(f)} \, {\rm trace}_g ({\mathcal L}_V
g)_{S(T{\mathcal F})} + \frac{g(V,V)}{2 V(f)^2} [\square f - 2
{\rm Hess}_f (\xi , N)] \] and (\ref{e:fol44}) is proved. Let us
look at the example (\ref{e:fol45}). The tangent bundle
$T({\mathcal F})$ is the span of
\[ \{ \partial_i - \sqrt{\frac{n-s}{s}} \, \partial_n \; , \;
\partial_j - \partial_n \, : 1 \leq i \leq s, \;\;\; s+1 \leq j
\leq n-1 \} . \] Moreover $\xi = - \sqrt{\frac{n-s}{s}} \,
\sum_{i=1}^s \partial_i + \sum_{j=s+1}^n \partial_j$ hence $\xi
\in T({\mathcal F})$ and $\{ X_a , \xi : 1 \leq a \leq n-2 \}$ are
linearly independent everywhere on $M$, where the $X_a$'s are
given by (\ref{e:fol46}). Therefore $S(T{\mathcal F})$ is indeed a
screen distribution and $S(T{\mathcal F})^\bot$ is the span of
\[ \{ \partial_s \; , \; - \sqrt{\frac{n-s}{s}} \,
\sum_{i=1}^{s-1} \partial_i + \sum_{j=s+1}^n \partial_j \} . \]
Consequently $V \in S(T{\mathcal F})^\bot$ where $V$ is given by
(\ref{e:fol47}). Next $\{ \xi , V \}$ are independent so that
(\ref{e:fol40}) holds. For any $C^\infty$ function $f: {\mathbb
R}^n_s \to {\mathbb R}$ one has ${\mathcal L}_\xi g = 2 {\rm
Hess}_f$. Finally a calculation shows that
\[ N = \frac{1}{2} \sqrt{\frac{s}{n-s}} \{ - \sum_{i=1}^{s-1}
\partial_i +  \partial_s + \sqrt{\frac{s}{n-s}} \sum_{j=s+1}^n
\partial_j \} \] hence (by (\ref{e:fol44})) $\kappa (N) = 0$.

\end{document}